\documentclass[11pt, a4paper]{amsart}
\usepackage[top=35truemm,bottom=35truemm,left=35truemm,right=35truemm]{geometry}
\usepackage{setspace} 
\usepackage{amsmath, amsfonts,amsthm,amssymb,mathrsfs}
\usepackage{enumerate}
\usepackage[pdftex]{graphicx}
\usepackage{ascmac}
\usepackage[arrow,matrix]{xy}
\usepackage{color}
\usepackage{array}
\usepackage{pifont}
\usepackage{longtable}
\usepackage{mathtools}

\theoremstyle{definition}
\newtheorem{theo}{Theorem}
\newtheorem*{theo*}{Theorem}
\newtheorem{defi}{Definition}

\newtheorem{prop}{Proposition}
\newtheorem{cor}{Corollary}
\newtheorem{lem}{Lemma}

\begin{document}

\title[ ]{On the growth rates of cofinite 3-dimensional hyperbolic Coxeter groups whose dihedral angles are of the form $\frac{\pi}{m}$ for $m=2,3,4,5,6$}
\author{Tomoshige Yukita}
\address{Department of Mathematics, School of Education, Waseda University, Nishi-Waseda 1-6-1, Shinjuku, Tokyo 169-8050, Japan}
\email{yshigetomo@suou.waseda.jp}
\date{}
\thanks{}

\begin{abstract}
We study arithmetic properties of the growth rates of cofinite 3-dimensional hyperbolic Coxeter groups whose dihedral angles are of the form $\frac{\pi}{m}$ for $m=2,3,4,5,6$ and show that the growth rates are always Perron numbers.
\end{abstract}

\maketitle

\setstretch{1.1}
\section{Introduction}
Let $\mathbb{H}^{n}$ denote the upper half-space model of hyperbolic $n$-space and $\overline{\mathbb{H}}^{n}$ its closure in $\mathbb{R}^n\cup{\{\infty\}}$. 
A convex polyhedron $P\subset{\overline{\mathbb{H}}^{n}}$ of finite volume is called a \textit{Coxeter polyhedron} if all of its dihedral angles are of the form $\frac{\pi}{m}$ for an integer $m\geq{2}$ or $m=\infty$.  
The set $S$ of reflections with respect to the facets of $P$ generates the discrete group $\Gamma$, called the \textit{hyperbolic Coxeter group},
and the pair $(\Gamma, S)$ is called the \textit{Coxeter system} associated with $P$. 
Then $P$ becomes a fundamental domain of $\Gamma$. 
If $P$ is compact (resp. non-compact), the hyperbolic Coxeter group $\Gamma$ is called \textit{cocompact} (resp. \textit{cofinite}). 
The \textit{growth series} $f_{S}(t)$ of $(\Gamma, S)$ is the formal power series $\sum_{k=0}^{\infty}a_{k}t^{k}$ where $a_{k}$ is the number of elements of $\Gamma$ whose word lengths with respect to $S$ are equal to $k$. 
Then $\tau_{\Gamma} :=\limsup_{k \rightarrow \infty} \sqrt[k]{a_k}$ is called the \textit{growth rate} of $(\Gamma, S)$. 
By means of Cauchy-Hadamard formula, $\tau_{\Gamma}$ is equal to the reciprocal of the radius of convergence $R$ of $f_{S}(t)$. 
It is known that the growth rate $\tau_{\Gamma}$ is a real algebraic integer bigger than 1 (\cite{dlH1}).
Certain classes of real algebraic integers show up in the study of the growth rates of cocompact or cofinite hyperbolic Coxeter group.
For two and three-dimensional cocompact hyperbolic Coxeter groups, Cannon-Wagreich and Parry showed that their growth rates are Salem numbers (\cite{CW}, \cite{Pa}), where a real algebraic integer $\tau>1$ is called a \textit{Salem number} if $\tau^{-1}$ is an algebraic conjugate of $\tau$ and all other algebraic conjugates lie on the unit circle.
Floyd also showed that the growth rates of two-dimensional cofinite hyperbolic Coxeter groups are Pisot-Vijayaraghavan numbers(\cite{F}), where a real algebraic integer $\tau>1$ is called a \textit{Pisot-Vijayaraghavan} number if all its algebraic conjugates are less than $1$ in absolute value.
In general, Kellerhals and Perren conjectured that the growth rates of cocompact hyperbolic Coxeter groups are Perron numbers (\cite{KP}), where a real algebraic integer $\tau>1$ is called a \textit{Perron number} if all its algebraic conjugates are less than $\tau$ in absolute value.
In the study of three-dimensional cofinite hyperbolic Coxeter groups, Komori and Umemoto showed that the growth rates of three-dimensional cofinite hyperbolic Coxeter simplex groups  are Perron numbers (\cite{KU}). 
Kellerhals and Nonaka showed that the growth rates of three-dimensional ideal hyperbolic Coxeter groups are Perron  numbers (\cite{N}); 
a hyperbolic Coxeter group is called ideal if all of the vertices of its fundamental polyhedron are located on the ideal boundary of hyperbolic space.
Komori and Yukita also showed the same result independently (\cite{KY}). 
In this paper, we consider the growth rates of cofinite Coxeter groups in hyperbolic 3-space and we shall prove the following theorem.
\begin{theo*}
Suppose $P$ is a non-compact hyperbolic Coxeter polyhedron whose dihedral angles are of the form $\frac{\pi}{m}$ for $m=2,3,4,5,6$. 
Then the growth rate of $P$ is a Perron number.
\end{theo*}

The organization of the present paper is as follows.
In section 2, we introduce useful formulas to calculate the growth series $f_{S}(t)$ and combinatorial identities with respect to a hyperbolic Coxeter polyhedron.
Then, we calculate the growth series $f_{S}(t)$ of a non-compact hyperbolic Coxeter polyhedron whose dihedral angles are of the form $\frac{\pi}{m}$ for $m=2,3,4,5,6$ in section 3.

\section{Preliminaries}
In this section,  we introduce some notations and useful identities to calculate the growth functions of hyperbolic Coxeter groups.

\begin{defi}{\rm(Coxeter system, Coxeter graph, growth rate)}

{\rm{(i)}} A {\textit {Coxeter system}} $(W,S)$ consists of a group $W$ and a set of generators $S\subset{W}$ with relations $(s_is_j)^{m_{ij}}$ for each $i , j$ , where $m_{ii}=1$ and $m_{ij}\geq{2}$ or $m_{ij}=\infty$ for $i\neq{j}$.
We call $W$ a {\textit {Coxeter group}}.
For any subset $I\subset{S}$, we define $W_I$ to be the subgroup of $W$ generated by $I$. $W_I$ is called the {\textit {Coxeter subgroup}} of $W$.

{\rm{(ii)}} The Coxeter graph of $(W,S)$ is constructed as follows: \\
the vertex set represents $S$. 
If $m_{ij}\geq{3}$ for $s_i\neq{s_j} \in{S}$, join the corresponding pair of vertices by an edge and label such an edge with $m_{ij}$.

{\rm{(iii)}} The {\textit{growth series}} $f_{S}(t)$ of $(W,S)$ is the formal power series $\sum_{k=0}^{\infty}a_kt^k$ where $a_k$ is the number of elements of $W$ whose word lengths with respect to $S$ are equal to $k$. Then $\tau =\limsup_{k \rightarrow \infty} \sqrt[k]{a_k}$ is called the {\textit{growth rate}} of $(W,S)$.

\end{defi}

A Coxeter group $W$ is \textit{irreducible} if the Coxeter graph of $(W,S)$ is connected.
In this paper, we are interested in Coxeter groups which acts on hyperbolic $n$-space as a discrete group. 

\begin{defi}{\rm (hyperbolic $n$-space)}

The upper half-space $\mathbb{H}^n:=\{x=(x_1,\cdots, x_n)\in{\mathbb{R}^n} \mid x_{n}>0\}$ with the metric $\frac{|dx|}{x_n}$ is called the \textit{upper half-space model of hyperbolic $n$-space}. 
$\overline{\mathbb{H}}^{n}$ (resp. $\partial{\mathbb{H}^n}$) denotes the closure (resp. the boundary) of $\mathbb{H}^n$ in $\mathbb{R}^n\cup{\{\infty\}}$. 
$\partial{\mathbb{H}^n}$ is called the \textit{ideal boundary of} $\mathbb{H}^n$.  
\end{defi}

\begin{defi}{\rm(hyperbolic polyhedron)} 

A subset $P\subset{\overline{\mathbb{H}}^n}$ is called a \textit{hyperbolic polyhedron} if $P$ can be written as the intersection of finitely many closed half spaces: $P=\cap{H_S}$ where $H_S$ is the closed domain of $\mathbb{H}^n$ bounded by the hyperplane $S$.
\end{defi}

Suppose that two hyperplanes $S$ and $T$ bounding the polyhedron $P$ satisfy $S\cap{T}\neq{\emptyset}$  in $\mathbb{H}^{n}$.
Then we define the {\textit {dihedral angle}} between $S$ and $T$ as follows: let us choose a point $x\in{S\cap{T}}$ and consider the outer-normal vectors $e_S, e_T\in{\mathbb{R}^n}$ of $S$ and $T$ with respect to $P$ starting from $x$. 
Then the dihedral angle between $S$ and $T$ is defined by the real number $\theta\in{[0, \pi)}$ satisfying $\cos{\theta}=-(e_S, e_T)$ where $(\cdot , \cdot)$ denote the Euclidean inner product of $\mathbb{R}^n$.

If $S\cap{T}=\emptyset$ in $\mathbb{H}^n$ and $\overline{S}\cap{\overline{T}}$ in $\overline{\mathbb{H}}^{n}$ is a point at the ideal boundary of $\mathbb{H}^n$, we define the dihedral angle between $S$ and $T$ is equal to zero.

\begin{defi}{\rm(hyperbolic Coxeter polyhedron)}
 
A hyperbolic polyhedron $P\subset{\overline{\mathbb{H}}^{n}}$ of finite volume is called a {\textrm{hyperbolic Coxeter polyhedron}} if all of its dihedral angles are of the form $\frac{\pi}{m}$ for an integer $m\geq{2}$ or $m=\infty$.
\end{defi}

Note that a hyperbolic polyhedron $P\subset{\overline{\mathbb{H}}^{n}}$ is of finite volume if and only if $P\cap{\partial{\mathbb{H}^n}}$ consists of vertices of $P$.
If $P\subset{\overline{\mathbb{H}}^{n}}$ is a hyperbolic Coxeter polyhedron, the set $S$ of reflections with respect to the facets of $P$ generates the discrete group $\Gamma$. We call $\Gamma$ the {\textit{$n$-dimensional hyperbolic Coxeter group}} associated with $P$. Moreover, if $P$ is compact (resp. non-compact), $\Gamma$ is called {\textit{cocompact}} (resp. {\textit{cofinite}}).

We recall Solomon's formula and Steinberg's formula which are very useful for calculating the growth series.

\begin{theo}{\rm (Solomon's formula)}\cite{So} \\
The growth series $f_S(t)$ of an irreducible finite Coxeter group $(\Gamma, S)$ can be written as 
$f_S(t)=[m_1+1,m_2+1, \cdots, m_k+1]$ where $[n] =1+t+ \cdots +t^{n-1}, [m,n]=[m][n]$,etc. and $\{m_1, m_2, \cdots, m_k \}$
is the set of exponents of $(\Gamma, S)$.
\end{theo}

Table 1 shows the exponents of irreducible finite Coxeter groups (see \cite{Hu} for details).

\begin{table}[h]
\begin{center}
\caption{Exponents}
\begin{tabular}{|c|c|c|}
\hline 
Coxeter group & Exponents & growth series \\ 
\hline
$A_n$ & $1,2,\cdots,n$ & $[2,3,\cdots,n+1]$ \\ 
\hline
$B_n$ & $1,3,\cdots,2n-1$ & $[2,4,\cdots,2n]$ \\ 
\hline
$D_n$ & $1,3,\cdots,2n-3,n-1$ & $[2,4,\cdots,2n-2][n]$ \\ 
\hline
$E_6$ & 1,4,5,7,8,11 & [2,5,6,8,9,12] \\ 
\hline
$E_7$ & 1,5,7,9,11,13,17 & [2,6,8,10,12,14,18] \\ 
\hline
$E_8$ & 1,7,11,13,17,19,23,29 & [2,8,12,14,18,20,24,30] \\ 
\hline
$F_4$ & 1,5,7,11 & [2,6,8,12] \\ 
\hline
$H_3$ & 1,5,9 & [2,6,10] \\ 
\hline
$H_4$ & 1,11,19,29 & [2,12,20,30] \\ 
\hline
$I_2(m)$ & 1,$m$-1 & [2,$m$] \\ 
\hline
\end{tabular}
\end{center}
\end{table}

\begin{theo}{\rm (Steinberg's formula)}\cite{St} \\
Let $(\Gamma, S)$ be a Coxeter group.
Let us denote the Coxeter subgroup of $(\Gamma, S)$ generated by the subset $T\subseteq S$ by $(\Gamma_T,T)$, and denote its growth series by $f_T(t)$.
Set $\mathcal{F}=\{T\subseteq S \;:\; \Gamma_T$ is finite $\}$. Then
$$
\frac{1}{f_S(t^{-1})}=\sum _{T \in \mathcal{F}} \frac{(-1)^{|T|}}{f_T(t)}.
$$
\end{theo}

By Theorem 1 and Theorem 2,  there exists a rational function $\frac{p(t)}{q(t)} (p,q\in{\mathbb{Z}[t]})$ such that $\frac{p(t)}{q(t)}$ is the analytic continuation of $f_{S}(t)$. 
The rational function $\frac{p(t)}{q(t)}$ is called the \textit{growth function} of $(\Gamma, S)$. 
The radius of convergence $R$ of the growth series $f_{S}(t)$ is equal to the real positive root of $q(t)$ which has the smallest absolute value of all roots of $q(t)$. 

Let $P$ be a hyperbolic Coxeter polyhedron in $\mathbb{H}^{3}$. 
The \textit {hyperbolic sphere} of $\mathbb{H}^{3}$ with center $a$ and radius $r>0$ is defined to be the set $S(a, r)=\{x\in{\mathbb{H}^{3}} \mid d_{H}(a,x)=r \}$ where $d_{H}$ denotes the hyperbolic metric on the upper half-space model of hyperbolic $n$-space. 
A subset $\Sigma$ of $\mathbb{H}^{3}$ is called a \textit {horosphere} of $\mathbb{H}^{3}$ based at a point $b$ of $\partial \mathbb{H}^{3}$ if $\Sigma$ is either a Euclidean hyperplane in $\mathbb{H}^{3}$ parallel to $E^2:=\{(x_1, x_2, 0) | x_1, x_2\in{\mathbb{R}}\}$ if $b=\infty$, or the intersection with $\mathbb{H}^{3}$ of a Euclidean sphere in $\overline{\mathbb{H}^{3}}$ tangent to $E^2$ at $b$ if $b\neq{\infty}$. 
When we restrict the hyperbolic metric of $\mathbb{H}^{3}$ to a hyperbolic sphere $S(a, r)$ (resp. a horosphere $\Sigma$ ), it becomes a model of 2-dimensional spherical geometry (resp. Euclidean geometry) .  

\begin{lem} ([10, Theorem 6.4.1 , 6.4.5])
Suppose $P$ is a hyperbolic Coxeter polyhedron in $\overline{\mathbb{H}^{3}}$ and $v$ is a vertex of $P$.
Let $S(v,r)$ (resp. $\Sigma$) be a hyperbolic sphere (resp. a horosphere) of $\overline{\mathbb{H}^{3}}$ with center $v$ such that $S(v,r)$ (resp. $\Sigma$) meets just the facets of $P$ incident to $v$ if $v$ is a finite vertex (resp. an ideal vertex).
Then the \textit{vertex link} $L(v):=P\cap{S(a, r)}$ (resp. $L(v):=P\cap{\Sigma}$) of $v$ in $P$ is a spherical (resp.  Euclidean) Coxeter polygon in the $S(v, r)$ (resp. $\Sigma$).
Moreover, if $F_{1}$ and $F_{2}$ are adjacent facets of $P$ incident to $v$, the hyperbolic dihedral angle between $F_{1}$ and $F_{2}$ is equal to the spherical (resp. Euclidean) interior angle between $F_{1}\cap{S(v, r)}$ and $F_{2}\cap{S(v, r)}$ (resp. $F_{1}\cap{\Sigma}$ and $F_{2}\cap{\Sigma}$).
\end{lem}

By Lemma 1, we obtain a necessary condition of dihedral angles of a hyperbolic Coxeter polyhedron $P$.

\begin{cor}
Suppose $P\subset{\overline{\mathbb{H}^{3}}}$ is a hyperbolic Coxeter polyhedron and $v$ is a vertex of $P$.
Let $F_{1},\cdots, F_{n}$ be the adjacent facets of $P$ incident to $v$ and $\frac{\pi}{a_{i}}$ (resp. $\frac{\pi}{a_n}$) be the hyperbolic dihedral angle between $F_{i}$ and $F_{i+1}$ (resp. $F_n$ and $F_1$).
Then, the number of facets of $P$ incident to $v$ is at most 4 and $a_1,\cdots, a_n$ satisfies the following. 
\begin{eqnarray}
a_{1}=a_{2}=a_{3}=a_{4}=2 & \qquad & \text{if $n=4$}. \\
\dfrac{1}{a_{1}}+\dfrac{1}{a_{2}}+\dfrac{1}{a_{3}}\geq{1}  & & \text{if $n=3$} .
\end{eqnarray}
\end{cor}

\proof By Lemma 1, the vertex link $L(v)$ is a spherical or Euclidean Coxeter $n$-sided polygon and the spherical or Euclidean interior angle is equal to the hyperbolic dihedral angle at the edge of $P$ incident to $v$.
By using Gauss-Bonnet theorem, 
\begin{center}
$\dfrac{\pi}{a_{1}}+\cdots+\dfrac{\pi}{a_{n}}\geq{(n-2)\pi}$.
\end{center}
Since $P$ is a Coxeter polyhedron, $2\leq{a_i}$ for all $i$. 
Then, this inequality implies that 
\begin{center}
$\dfrac{n}{2}\geq{n-2}$ .
\end{center}
Thus, we conclude that $n\leq{4}$ and obtain the identity (1) and inequality (2).
\qed

Note that a vertex $v$ of $P$ is an ideal vertex if and only if $a_{1}=a_{2}=a_{3}=a_{4}=2$ or $\dfrac{1}{a_{1}}+\dfrac{1}{a_{2}}+\dfrac{1}{a_{3}}=1$, and we call an ideal vertex a \textit{cusp}. 
\begin{itemize}
 \item If a vertex $v$ of $P$ satisfies the identity $(1)$, we call $v$ a \textit{cusp of type $(2,2,2,2)$}. 
 \item If a vertex $v$ of $P$ satisfies the inequality $(2)$, we call $v$ a \textit{vertex of type $(a_1, a_2, a_3)$} .
 \item $v_{2,2,2,2}$ denotes the number of cusps of type $(2,2,2,2)$. 
 \item $v_{a_1, a_2, a_3}$ denotes the number of vertices of type $(a_1, a_2, a_3)$. 
 \item $v, e, f$ denotes the number of vertices, edges and facets of $P$. 
 \item If an edge $e$ of $P$ has dihedral angle $\frac{\pi}{m}$,  we call it \textit{$\frac{\pi}{m}$-edge}. 
 \item $e_m$ denotes the number of $\frac{\pi}{m}$-edges. 
\end{itemize}

\begin{lem}
Let $P\subset{\overline{\mathbb{H}^3}}$ be a non-compact hyperbolic Coxeter polyhedron and $f_{P}(t)$ be the growth function of $P$. Then the following identities and inequality hold. 
\begin{eqnarray} 
 & &v-e+f= 2 .\\
 & & v=v_{2,2,2,2}+\sum_{n\geq{2}}^{}v_{2,2,n}+v_{2,3,3}+v_{2,3,4}+v_{2,3,5}+v_{2,3,6}+v_{2,4,4}+v_{3,3,3} .\\
 & & e=\sum_{n\geq2}^{}e_{n} .\\ 
 & &2e_{2} = 4v_{2,2,2,2}+3v_{2,2,2}+2\sum_{n\geq{3}}^{}v_{2,2,n}+v_{2,3,3}+v_{2,3,4}+v_{2,3,5}+v_{2,3,6}+v_{2,4,4} .\\
 & &2e_{3} = 3v_{3,3,3}+2v_{2,3,3}+v_{2,2,3}+v_{2,3,4}+v_{2,3,5}+v_{2,3,6} .\\
 & &2e_{4} = 2v_{2,4,4}+v_{2,2,4}+v_{2,3,4} .\\
 & &2e_{5} = v_{2,2,5}+v_{2,3,5} .\\
 & &2e_{6} = v_{2,2,6}+v_{2,3,6} .\\
 & &2e_{n} = v_{2,2,n}  \qquad \text{$n\geq{7}$} .\\
 & & v_{2,2,2,2}+v_{2,3,6}+v_{2,4,4}+v_{3,3,3}\geq{1} .
\end{eqnarray}
\end{lem}

\proof The identity $(3)$ is Euler's polyhedral formula. By the definitions of $v, e$ and $e_n$, the identities $(4)$ and $(5)$ hold for $P$. 
Non-compactness of $P$ implies the inequality (12).
The other identities are obtained by counting the number of edges which is adjacent to each vertex. 
For example, the identity (8) is obtained as follows. 
Any $\frac{\pi}{4}$-edge has strictly two vertices of type $(2,2,4)$ or $(2,3,4)$ or $(2,4,4)$ (see Fig 1). 
On the other hand, any vertex of type $(2,2,4)$ or $(2,3,4)$ has one $\frac{\pi}{4}$-edge and any vertex of type $(2,4,4)$ has two $\frac{\pi}{4}$-edges. 
Hence, we obtain the identity (8).

\begin{figure}[htbp]
\begin{center}
 \includegraphics [width=165pt, clip]{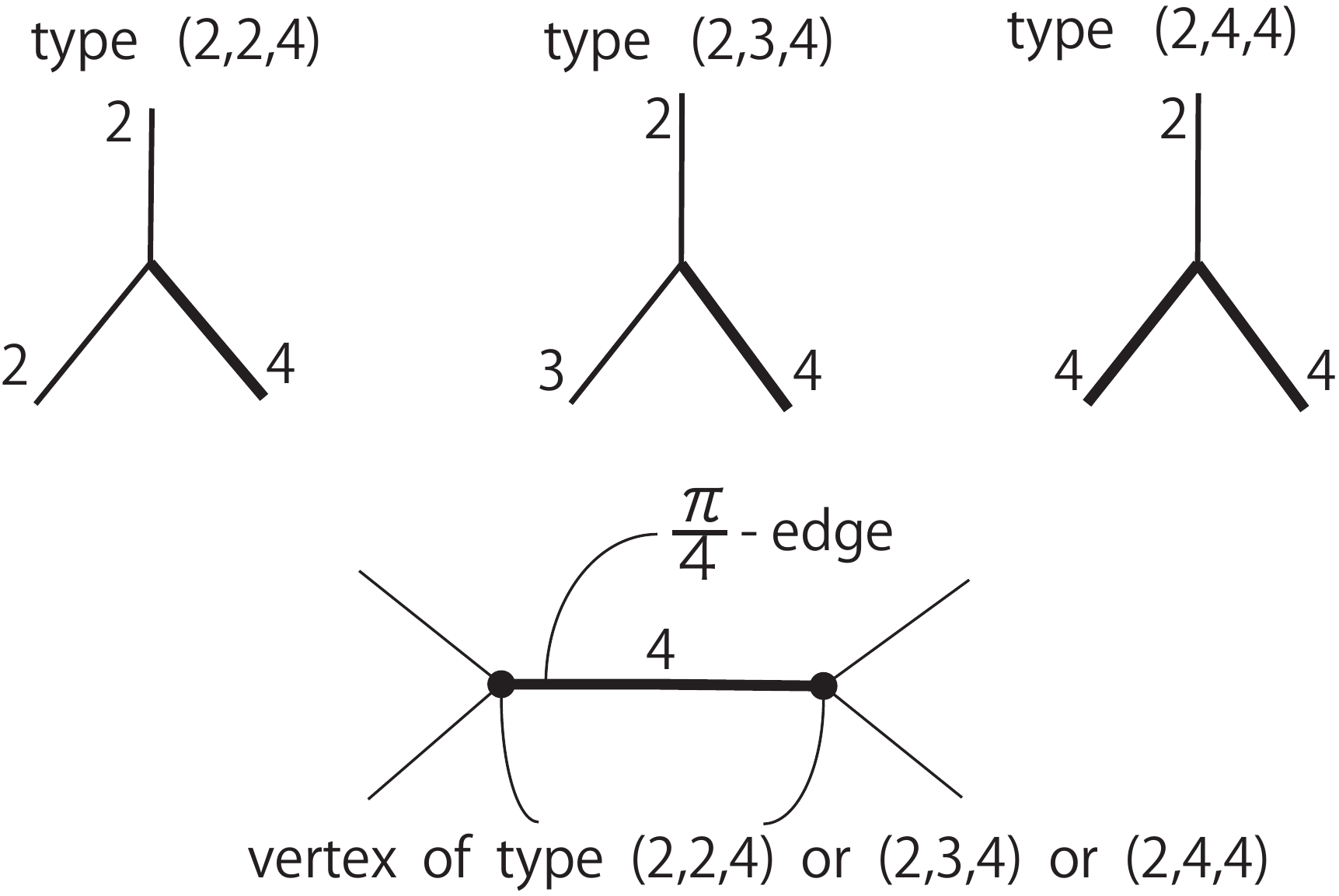}
\end{center}
\caption{}
\end{figure} 
\qed

We use the following result to show that the growth rates of 3-dimensional cofinite hyperbolic Coxeter groups are Perron numbers.

\begin{prop}\rm{(\cite{KU}, Lemma1)}\\
\label{prop:KU}
 Consider the following polynomial of degree $n \geq 2$
 $$
 g(t)=\sum _{k=1}^n a_k t^k -1,
 $$
 where  $a_k$ is a nonnegative integer.
 We also assume that the greatest common divisor of $\{k \in \mathbb{N} \; |\; a_k \neq 0 \}$ is $1$.
 Then there is a real number $r_0$,  $0<r_0<1$ which is the unique zero of $g(t)$ having the smallest absolute value of all zeros of $g(t)$.
\end{prop}

\section{non-compact Coxeter polytopes whose dihedral angles are of the form $\frac{\pi}{m}$ for $m=2,3,4,5,6$ }
In this section, we calculate the growth function $f_{P}(t)$ of a hyperbolic Coxeter group associated with a non-compact hyperbolic Coxeter polyhedron $P$ whose dihedral angles are of the form $\frac{\pi}{m}$ for $m=2,3,4,5,6$ and prove that the growth rate of $P$ is a Perron number. 
In the sequel, we call $f_{P}(t)$ the growth function of $P$.
Note that Komori and Umemoto proved that the growth rates of 3-dimensional cofinite hyperbolic Coxeter simplex groups are Perron numbers (\cite{KU}). For this reason, we consider non-compact hyperbolic Coxeter polyhedra with at least 5-facets.

\subsection{The growth rates of non-compact Coxeter polyhedra whose dihedral angles are $\frac{\pi}{2}$}

Let $P$ be a non-compact hyperbolic Coxeter polyhedron whose dihedral angles are $\frac{\pi}{2}$.
By means of Steinberg's formula, we can calculate the growth function $f_{P}(t)$ of $P$ as

\[
\dfrac{1}{f_{P}(t^{-1})}=1-\dfrac{f}{[2]}+\dfrac{e_{2}}{[2,2]}-\dfrac{v_{2,2,2}}{[2,2,2]}.
\]

By using the identities (3), (4), (5) and (6), it can be rewritten as

\begin{eqnarray*}
\dfrac{1}{f_{P}(t)} &=& 1-\dfrac{ft}{[2]}+\dfrac{e_{2}t^{2}}{[2,2]}-\dfrac{v_{2,2,2}t^{3}}{[2,2,2]}\\
                               &=& \dfrac{1}{[2,2,2]}\Bigl\{ (v_{2,2,2,2}-1)t^2+(v_{2,2,2,2}+\dfrac{1}{2}v_{2,2,2}-2)t-1\Bigr\}\\
                               &=&\dfrac{(t-1)}{[2,2,2]}\Bigl\{ (v_{2,2,2,2}-1)t^2+(f-4)t-1\Bigr\}\\
                               &=& \dfrac{(t-1)}{[2,2,2]}H_{2}(t)\\
\end{eqnarray*}
where we put
\begin{eqnarray*}
H_{2}(t) &=& (v_{2,2,2,2}-1)t^2+(f-4)t-1.
\end{eqnarray*}

\begin{prop}
All the coefficients of $H_{2}(t)$ is nonnegative except its constant term. Moreover, the growth rate of $P$ is a Perron number.
\end{prop}

\proof
Put $a_{2}=v_{2,2,2,2}-1$ and $a_{1}=f-4$.
$f\geq{5}$ and the inequality (12) imply that $a_1, a_2\geq{0}$.
If $v_{2,2,2,2}=1$, $H_{2}(t)$ has only one real positive root which is less than 1. Hence, we conclude that the growth rate of $P$ is a Perron number.
If $v_{2,2,2,2}\geq{2}$, by using Proposition 1, we conclude that the growth rate of $P$ is a Perron number.
\qed

\subsection{The growth rates of non-compact Coxeter polyhedra whose dihedral angles are $\frac{\pi}{2}$ and $\frac{\pi}{3}$}
Let $P$ be a non-compact hyperbolic Coxeter polyhedron whose dihedral angles are $\frac{\pi}{2}$ and $\frac{\pi}{3}$.
We also assume that $P$ has at least one $\frac{\pi}{3}$-edge.
By means of Steinberg's formula and using the identity , we calculate the growth function of $P$ as

\begin{eqnarray*}
\dfrac{1}{f_{P}(t)} &=& \dfrac{(t-1)}{[2,2,3,4]}H_{2,3}(t)\\
\end{eqnarray*}
where
\begin{eqnarray*}
H_{2,3}(t) &=& (v_{2,2,2,2}+v_{3,3,3}-1)t^6\\
                  &+& (2v_{2,2,2,2}+\frac{1}{2}v_{2,2,2}+\frac{1}{2}v_{2,2,3}+\frac{1}{2}v_{2,3,3}+\frac{5}{2}v_{3,3,3}-3)t^5\\
                  &+& (3v_{2,2,2,2}+\frac{1}{2}v_{2,2,2}+v_{2,2,3}+\frac{3}{2}v_{2,3,3}+3v_{3,3,3}-5)t^4\\
                  &+& (3v_{2,2,2,2}+v_{2,2,2}+v_{2,2,3}+2v_{2,3,3}+3v_{3,3,3}-6)t^3\\
                  &+& (2v_{2,2,2,2}+\frac{1}{2}v_{2,2,2}+v_{2,2,3}+\frac{3}{2}v_{2,3,3}+2v_{3,3,3}-5)t^2\\
                  &+& (v_{2,2,2,2}+\frac{1}{2}v_{2,2,2}+\frac{1}{2}v_{2,2,3}+\frac{1}{2}v_{2,3,3}+\frac{1}{2}v_{3,3,3}-3)t-1.\\
\end{eqnarray*}

By using the identities (3), (4), (5), (6) and (7),  $H_{2,3}(t)$ can be rewritten as

\begin{eqnarray*}
H_{2,3}(t) &=& (v_{2,2,2,2}+v_{3,3,3}-1)t^6\\
                  &+& (v_{2,2,2,2}+2v_{3,3,3}+f-5)t^5\\
                  &+& (2v_{2,2,2,2}+\frac{1}{2}v_{2,2,3}+v_{2,3,3}+\frac{5}{2}v_{3,3,3}+f-7)t^4\\
                  &+& (v_{2,2,2,2}+v_{2,3,3}+2v_{3,3,3}+2f-10)t^3\\
                  &+& (v_{2,2,2,2}+\frac{1}{2}v_{2,2,3}+v_{2,3,3}+\frac{3}{2}v_{3,3,3}+f-7)t^2\\
                  &+& (f-5)t-1.\\
\end{eqnarray*}

\begin{prop}
All the coefficients of $H_{2,3}(t)$ except its constant term are nonnegative and the growth rate of $P$ is a Perron number.
\end{prop}
\proof
We denote the coefficient of $k$-th term of $H_{2,3}(t)$ by $a_{k}$, that is,  $H_{2,3}(t)=\sum_{k=1}^{6}a_{k}t^k-1$.

By the inequality (12) with respect to $P$ and $f\geq{5}$, we obtain the following inequalities. 
\begin{eqnarray*}
a_{6}=v_{2,2,2,2}+v_{3,3,3}-1 &\geq& 0. \\
a_{5}=v_{2,2,2,2}+2v_{3,3,3}+f-5>f-5 &\geq&0. \\
a_{3}=v_{2,2,2,2}+v_{2,3,3}+2v_{3,3,3}+2f-10>2f-10 &\geq&0. \\
a_{1}=f-5 &\geq&0.  
\end{eqnarray*}
These inequalities mean that $a_{1}, a_{3}, a_{5}$ and $a_{6}$ are nonnegative. Moreover, $a_{3}$ and $a_{5}$ are positive.
\begin{eqnarray*}
a_{4}=2v_{2,2,2,2}+\frac{1}{2}v_{2,2,3}+v_{2,3,3}+\frac{5}{2}v_{3,3,3}+f-7&\geq&2(v_{2,2,2,2}+v_{3,3,3})+f-7 \\
                                                                                                                                   &\geq& 2\cdot1+f-7 \geq{0}
\end{eqnarray*}
implies that $a_{4}$ is nonnegative.

Since $P$ has at least one $\frac{\pi}{3}$-edge and $\frac{\pi}{3}$-edge has two vertices of type $(2,2,3)$ or $(2,3,3)$ or $(3,3,3)$ (see Fig 2), we obtain the following inequality.
\begin{center}
$v_{2,2,3}+v_{2,3,3}+v_{3,3,3}\geq{2}$.
\end{center}

\begin{figure}[htbp]
\begin{center}
 \includegraphics [width=200pt, clip]{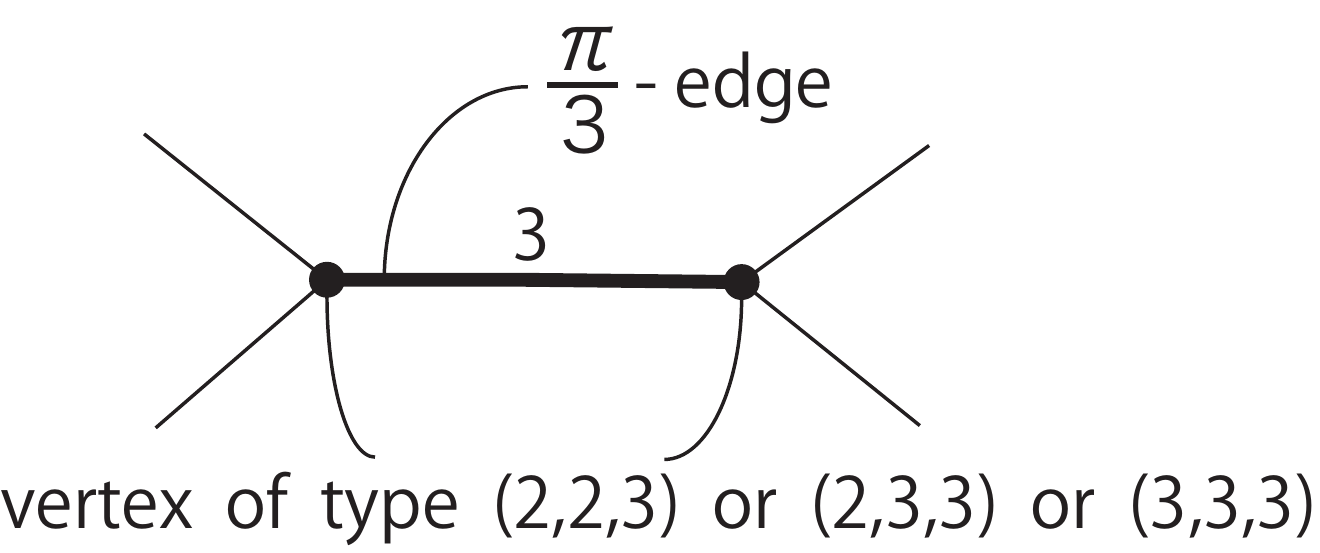}
\end{center}
\caption{}
\end{figure} 

Thus,
\begin{eqnarray*}
a_{2} &=& v_{2,2,2,2}+\frac{1}{2}v_{2,2,3}+v_{2,3,3}+\frac{3}{2}v_{3,3,3}+f-7\\
                                 &=& (v_{2,2,2,2}+v_{3,3,3})+\frac{1}{2}v_{2,2,3}+v_{2,3,3}+\frac{1}{2}v_{3,3,3}+f-7\\
                                 &\geq& (v_{2,2,2,2}+v_{3,3,3})+\frac{1}{2}(v_{2,2,3}+v_{2,3,3}+v_{3,3,3})+f-7\\
                                 &\geq& 1+\frac{1}{2}\cdot2+f-7\\
                                 &\geq& 0.
\end{eqnarray*}

Since $a_3$ and $a_5$ are positive and all the coefficients of $H_{2,3}(t)$ except its constant term are nonnegative,
we can apply Proposition 1 to $H_{2,3}(t)$; 
therefore, the growth rate of $P$ is a Perron number.  
\qed

\subsection{The growth rates of non-compact Coxeter polyhedra whose dihedral angles are $\frac{\pi}{2},\frac{\pi}{3}$ and $\frac{\pi}{6}$}
Let $P$ be a non-compact hyperbolic Coxeter polyhedron whose dihedral angels are $\frac{\pi}{2},\frac{\pi}{3}$ and $\frac{\pi}{6}$.
We also assume that $P$ has at least one $\frac{\pi}{6}$-edge. 
 By means of Steinberg's formula and using the identities (3), (4), (5), (6), (7) and (10), $f_{P}(t)$ can be written as
\begin{eqnarray*}
\dfrac{1}{f_{P}(t)} &=& \dfrac{(t-1)}{[2,2,4,6]}H_{2,3,6}(t)
\end{eqnarray*}
where
\begin{eqnarray*}
H_{2,3,6}(t) &=& (v_{2,2,2,2}+v_{2,3,6}+v_{3,3,3}-1)t^9\\
                      &+&(v_{2,2,2,2}+2v_{2,3,6}+2v_{3,3,3}+f-5)t^8\\
                      &+&(2v_{2,2,2,2}+\frac{1}{2}v_{2,2,3}+\frac{1}{2}v_{2,2,6}+v_{2,3,3}+3v_{2,3,6}+\frac{5}{2}v_{3,3,3}+f-7)t^7\\
                      &+&(2v_{2,2,2,2}+\frac{1}{2}v_{2,2,6}+v_{2,3,3}+\frac{7}{2}v_{2,3,6}+3v_{3,3,3}+2f-11)t^6\\
                      &+&(2v_{2,2,2,2}+\frac{1}{2}v_{2,2,3}+v_{2,2,6}+v_{2,3,3}+\frac{7}{2}v_{2,3,6}+\frac{7}{2}v_{3,3,3}+2f-12)t^5\\
                      &+&(2v_{2,2,2,2}+\frac{1}{2}v_{2,2,3}+v_{2,2,6}+v_{2,3,3}+\frac{5}{2}v_{2,3,6}+\frac{5}{2}v_{3,3,3}+2f-12)t^4\\
                      &+&(v_{2,2,2,2}+\frac{1}{2}v_{2,2,6}+v_{2,3,3}+\frac{3}{2}v_{2,3,6}+2v_{3,3,3}+2f-11)t^3\\
                      &+&(v_{2,2,2,2}+\frac{1}{2}v_{2,2,3}+\frac{1}{2}v_{2,2,6}+v_{2,3,3}+v_{2,3,6}+\frac{3}{2}v_{3,3,3}+f-7)t^2\\
                      &+&(f-5)t-1.
\end{eqnarray*}

\begin{prop}
All the coefficients of $H_{2,3,6}(t)$ except its constant term are non-negative and the growth rate is a Perron number.
\end{prop}

\proof
We denote the coefficient of $k$-th term of $H_{2,3,6}(t)$ by $a_{k}$, that is, $H_{2,3,6}(t)=\sum_{k=1}^{9}a_{k}t^k-1$.

The assumptions $f\geq{5}$ and that $P$ is non-compact imply
\[
a_{k}\geq{0} \ (k\neq{2}) .
\]
Since $P$ has at least one $\frac{\pi}{6}$-edge, by the identity (10), we obtain the following inequality.
\[
v_{2,2,6}+v_{2,3,6}=2e_{6}\geq{2}.
\]
For example, $a_7>0$ is obtained as follows.
\begin{eqnarray*}
a_{7} &=& 2v_{2,2,2,2}+\frac{1}{2}v_{2,2,3}+\frac{1}{2}v_{2,2,6}+v_{2,3,3}+3v_{2,3,6}+\frac{5}{2}v_{3,3,3}+f-7\\
                                 &=& 2(v_{2,2,2,2}+v_{2,3,6}+v_{3,3,3})+\frac{1}{2}v_{2,2,3}+v_{2,3,3}+\Bigl(\frac{1}{2}v_{2,2,6}+v_{2,3,6}\Bigr)+\frac{1}{2}v_{3,3,3}+f-7\\
                                 &\geq& 2\cdot1+1+f-7 \\
                                 &>& 0.
\end{eqnarray*}

For $a_2$, we obtain
\begin{eqnarray*}
a_2 &=& v_{2,2,2,2}+\frac{1}{2}v_{2,2,3}+\frac{1}{2}v_{2,2,6}+v_{2,3,3}+v_{2,3,6}+\frac{3}{2}v_{3,3,3}+f-7\\
        &=& (v_{2,2,2,2}+v_{2,3,6}+v_{3,3,3})+\frac{1}{2}v_{2,2,3}+\frac{1}{2}v_{2,2,6}+v_{2,3,3}+\frac{1}{2}v_{3,3,3}+f-7\\
        &\geq& (v_{2,2,2,2}+v_{2,3,6}+v_{3,3,3})+\frac{1}{2}v_{2,2,3}+\frac{1}{2}v_{2,2,6}+v_{2,3,3}+\frac{1}{2}v_{3,3,3}-2.
\end{eqnarray*}
Define $\tilde{a_{2}}$ to be the last expression in the inequality above.
The above argument can be split into the following three cases.
\begin{itemize}
 \item First, $P$ has at least two cusps.
 \item Second, $P$ has only one cusp of type $(2,2,2,2)$ or $(3,3,3)$.
 \item Finally, $P$ has only one cusp of type $(2,3,6)$.
\end{itemize}

Let us first consider the case that $P$ has at least two cusps. 
In this case, we obtain 
\begin{eqnarray*}
\tilde{a_2} &=& (v_{2,2,2,2}+v_{2,3,6}+v_{3,3,3})+\frac{1}{2}v_{2,2,3}+\frac{1}{2}v_{2,2,6}+v_{2,3,3}+\frac{1}{2}v_{3,3,3}-2 \\ 
                   &\geq& 2+\frac{1}{2}v_{2,2,3}+\frac{1}{2}v_{2,2,6}+v_{2,3,3}+\frac{1}{2}v_{3,3,3}-2\geq{0} .
\end{eqnarray*}

Next, we consider the case that $P$ has only one cusp of type $(2,2,2,2)$ or $(3,3,3)$.
All $\frac{\pi}{6}$-edges have two vertices of type $(2,2,6)$ or $(2,3,6)$.
 However, vertices of type $(2,3,6)$ are cusps and $P$ has only one cusp of type $(2,2,2,2)$ or $(3,3,3)$ in this case.
 Hence, all $\frac{\pi}{6}$-edges have two vertices of type $(2,2,6)$ (see Fig 3).
\begin{figure}[htbp]
\begin{center}
 \includegraphics [width=250pt, clip]{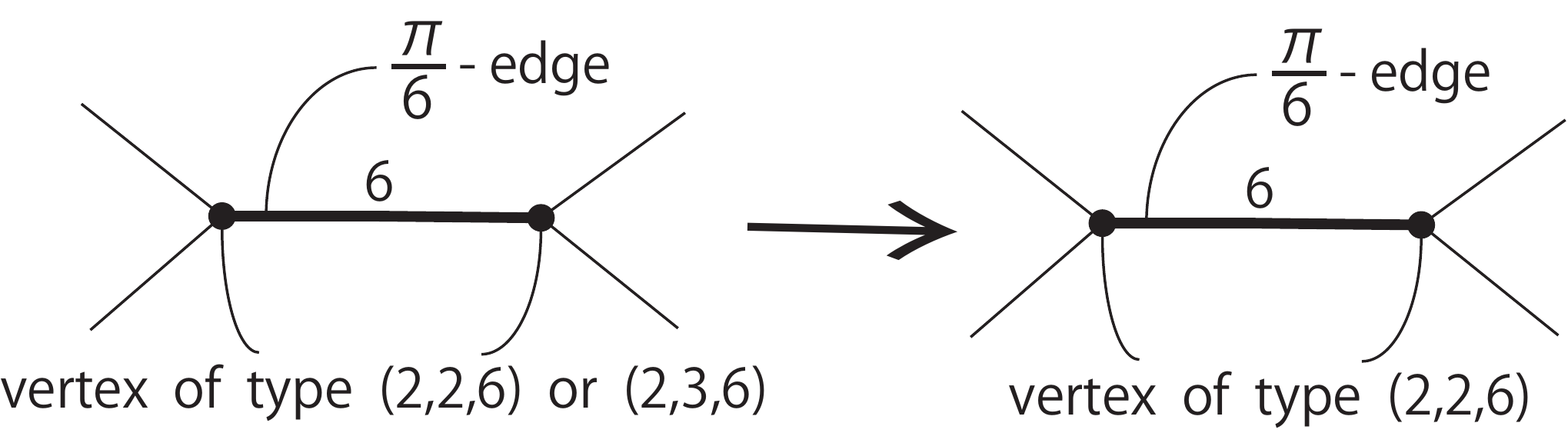}
\end{center}
\caption{}
\end{figure} 

 Since $P$ has at least one $\frac{\pi}{6}$-edge, we obtain the following inequality.
\[
v_{2,2,6}=2e_{6}\geq{2}.
\]
This inequality implies that 
\begin{eqnarray*}
\tilde{a_2} &=& (v_{2,2,2,2}+v_{2,3,6}+v_{3,3,3})+\frac{1}{2}v_{2,2,3}+\frac{1}{2}v_{2,2,6}+v_{2,3,3}+\frac{1}{2}v_{3,3,3}-2 \\ 
                   &\geq& 1+\frac{1}{2}v_{2,2,3}+\frac{1}{2}\cdot{2}+v_{2,3,3}+\frac{1}{2}v_{3,3,3}-2\geq{0} .
\end{eqnarray*}

Finally, we consider the case that $P$ has only one cusp of type $(2,3,6)$.
A cusp of type $(2,3,6)$ is shared by a $\frac{\pi}{2}$-edge, $\frac{\pi}{3}$-edge and $\frac{\pi}{6}$-edge.
Hence, $P$ has at least one $\frac{\pi}{3}$-edge.
By the identity (7), we obtain
\[
2v_{2,3,3}+v_{2,2,3}\geq{1} \\
\Leftrightarrow v_{2,3,3}+\frac{1}{2}v_{2,2,3}\geq{\frac{1}{2}}.
\]
Moreover, since $P$ has at least one $\frac{\pi}{6}$-edge and only one cusp of type $(2,3,6)$, we see that $P$ has at least one vertex of type $(2,2,6)$.
This means that $v_{2,2,6}\geq{1}$.
These facts show that 
\begin{eqnarray*}
\tilde{a_2} &=& (v_{2,2,2,2}+v_{2,3,6}+v_{3,3,3})+\frac{1}{2}v_{2,2,3}+\frac{1}{2}v_{2,2,6}+v_{2,3,3}+\frac{1}{2}v_{3,3,3}-2 \\ 
                   &\geq& 1+\frac{1}{2}+\frac{1}{2}+\frac{1}{2}v_{3,3,3}-2\geq{0} .
\end{eqnarray*}
We conclude that $\tilde{a_2}\geq{0}$; therefore all the coefficients of $H_{2,3,6}(t)$ except its constant term are nonnegative. Moreover, $a_7$ and $a_8$ are positive, so that we can apply Proposition 1 to $H_{2,3,6}(t)$.
Thus the growth rate of $P$ is a Perron number.
\qed

\subsection{The growth rates of non-compact Coxeter polyhedra whose dihedral angles are $\frac{\pi}{2},\frac{\pi}{3},\frac{\pi}{5}$ and $\frac{\pi}{6}$}
Let $P$ be a non-compact hyperbolic Coxeter polyhedron whose dihedral angels are $\frac{\pi}{2},\frac{\pi}{3},\frac{\pi}{5}$ and $\frac{\pi}{6}$.
We also assume that $P$ has at least one $\frac{\pi}{5}$-edge.
 By means of Steinberg's formula and using identity (3), (4), (5), (6), (7), (9) and (10), $f_{P}(t)$ can be written as
\begin{eqnarray*}
\dfrac{1}{f_{P}(t)} &=& \dfrac{(t-1)}{[2,4,6,10]}H_{2,3,5,6}(t)\\
\end{eqnarray*}
where
\begin{eqnarray*}
H_{2,3,5,6}(t) &=& (v_{2,2,2,2}+v_{2,3,6}+v_{3,3,3}-1)t^{17}\\
                      &+&(v_{2,2,2,2}+2v_{2,3,6}+2v_{3,3,3}+f-5)t^{16}\\
                      &+&(3v_{2,2,2,2}+\frac{1}{2}v_{2,2,3}+\frac{1}{2}v_{2,2,5}+\frac{1}{2}v_{2,2,6}+v_{2,3,3}+v_{2,3,5}+4v_{2,3,6}+\frac{7}{2}v_{3,3,3}+f-8)t^{15}\\
                      &+&(3v_{2,2,2,2}+\frac{1}{2}v_{2,2,5}+\frac{1}{2}v_{2,2,6}+v_{2,3,3}+\frac{3}{2}v_{2,3,5}+\frac{11}{2}v_{2,3,6}+5v_{3,3,3}+3f-16)t^{14}\\
                      &+& (5v_{2,2,2,2}+v_{2,2,3}+\frac{3}{2}v_{2,2,5}+\frac{3}{2}v_{2,2,6}+2v_{2,3,3}+\frac{7}{2}v_{2,3,5}+\frac{15}{2}v_{2,3,6}+7v_{3,3,3}+3f-20)t^{13}\\
                  &+&(5v_{2,2,2,2}+\frac{1}{2}v_{2,2,3}+v_{2,2,5}+\frac{3}{2}v_{2,2,6}+2v_{2,3,3}+\frac{9}{2}v_{2,3,5}+8v_{2,3,6}+\frac{15}{2}v_{3,3,3}+5f-28)t^{12}\\
         &+&(6v_{2,2,2,2}+v_{2,2,3}+2v_{2,2,5}+2v_{2,2,6}+3v_{2,3,3}+6v_{2,3,5}+9v_{2,3,6}+9v_{3,3,3}+5f-31)t^{11}\\
                 &+&(6v_{2,2,2,2}+v_{2,2,3}+\frac{3}{2}v_{2,2,5}+2v_{2,2,6}+3v_{2,3,3}+\frac{13}{2}v_{2,3,5}+9v_{2,3,6}+9v_{3,3,3}+6f-35)t^{10}\\
                  & +&(6v_{2,2,2,2}+v_{2,2,3}+2v_{2,2,5}+2v_{2,2,6}+3v_{2,3,3}+7v_{2,3,5}+9v_{2,3,6}+9v_{3,3,3}+6f-36)t^9\\
                  &+&(6v_{2,2,2,2}+v_{2,2,3}+2v_{2,2,5}+2v_{2,2,6}+3v_{2,3,3}+7v_{2,3,5}+9v_{2,3,6}+9v_{3,3,3}+6f-36)t^8\\
                  &+&(5v_{2,2,2,2}+v_{2,2,3}+\frac{3}{2}v_{2,2,5}+2v_{2,2,6}+3v_{2,3,3}+\frac{13}{2}v_{2,3,5}+8v_{2,3,6}+8v_{3,3,3}+6f-35)t^7\\
                  &+&(5v_{2,2,2,2}+v_{2,2,3}+2v_{2,2,5}+2v_{2,2,6}+3v_{2,3,3}+6v_{2,3,5}+7v_{2,3,6}+7v_{3,3,3}+5f-31)t^6\\
                  & +&(3v_{2,2,2,2}+\frac{1}{2}v_{2,2,3}+v_{2,2,5}+\frac{3}{2}v_{2,2,6}+2v_{2,3,3}+\frac{9}{2}v_{2,3,5}+5v_{2,3,6}+\frac{11}{2}v_{3,3,3}+5f-28)t^5\\
                  & +&(3v_{2,2,2,2}+v_{2,2,3}+\frac{3}{2}v_{2,2,5}+\frac{3}{2}v_{2,2,6}+2v_{2,3,3}+\frac{7}{2}v_{2,3,5}+\frac{7}{2}v_{2,3,6}+4v_{3,3,3}+3f-20)t^4\\
                   &+&(v_{2,2,2,2}+\frac{1}{2}v_{2,2,5}+\frac{1}{2}v_{2,2,6}+v_{2,3,3}+\frac{3}{2}v_{2,3,5}+\frac{3}{2}v_{2,3,6}+2v_{3,3,3}+3f-16)t^3\\
                  &+&(v_{2,2,2,2}+\frac{1}{2}v_{2,2,3}+\frac{1}{2}v_{2,2,5}+\frac{1}{2}v_{2,2,6}+v_{2,3,3}+v_{2,3,5}+v_{2,3,6}+\frac{3}{2}v_{3,3,3}+f-8)t^2\\
                   &+&(f-5)t-1.
\end{eqnarray*}

\begin{prop}
Suppose $P$ has at least one $\frac{\pi}{3}$-edge or $\frac{\pi}{6}$-edge. 
Then, 
all the coefficients of $H_{2,3,5,6}(t)$ except its constant term are nonnegative and the growth rate of $P$ is a Perron number.
\end{prop}

\proof
We denote the coefficient of $k$-th term of $H_{2,3,5,6}(t)$ by $a_{k}$, that is, $H_{2,3,5,6}(t)=\sum_{k=1}^{17}a_{k}t^k-1$.
 The assumptions $f\geq{5}$ and that $P$ is non-compact imply
\[
a_{k}\geq{0} \ (k\neq{2}) .
\]
We can see that $a_k$ $(k\neq{2})$ is nonnegative without the assumption that $P$ has at least one $\frac{\pi}{3}$-edge or $\frac{\pi}{6}$-edge. 
For example, $a_{15}\geq{0}$ is obtained as follows. By using the inequalities (12) and $f\geq{5}$, 
\begin{eqnarray*}
a_{15} &=& 3v_{2,2,2,2}+\frac{1}{2}v_{2,2,3}+\frac{1}{2}v_{2,2,5}+\frac{1}{2}v_{2,2,6}+v_{2,3,3}+v_{2,3,5}+4v_{2,3,6}+\frac{7}{2}v_{3,3,3}+f-8 \\
             &=& 3(v_{2,2,2,2}+v_{2,3,6}+v_{3,3,3})+\frac{1}{2}v_{2,2,3}+\frac{1}{2}v_{2,2,5}+\frac{1}{2}v_{2,2,6}+v_{2,3,3}+v_{2,3,5}+v_{2,3,6}+\frac{1}{2}v_{3,3,3}+f-8 \\
             &\geq& 3(v_{2,2,2,2}+v_{2,3,6}+v_{3,3,3})+\frac{1}{2}v_{2,2,3}+\frac{1}{2}v_{2,2,5}+\frac{1}{2}v_{2,2,6}+v_{2,3,3}+v_{2,3,5}+v_{2,3,6}+\frac{1}{2}v_{3,3,3}-3\\
             &\geq& \frac{1}{2}v_{2,2,3}+\frac{1}{2}v_{2,2,5}+\frac{1}{2}v_{2,2,6}+v_{2,3,3}+v_{2,3,5}+v_{2,3,6}+\frac{1}{2}v_{3,3,3} \geq{0}.
\end{eqnarray*}
Moreover, $e_5\neq{0}$ implies $a_{15}>0$. 
For $a_2$, we obtain 
\begin{eqnarray*}
a_{2} &=& v_{2,2,2,2}+\frac{1}{2}v_{2,2,3}+\frac{1}{2}v_{2,2,5}+\frac{1}{2}v_{2,2,6}+v_{2,3,3}+v_{2,3,5}+v_{2,3,6}+\frac{3}{2}v_{3,3,3}+f-8 \\
           &=& (v_{2,2,2,2}+v_{2,3,6}+v_{3,3,3})+\frac{1}{2}v_{2,2,3}+\frac{1}{2}v_{2,2,5}+\frac{1}{2}v_{2,2,6}+v_{2,3,3}+v_{2,3,5}+\frac{1}{2}v_{3,3,3}+f-8 \\ 
           &\geq& (v_{2,2,2,2}+v_{2,3,6}+v_{3,3,3})+\frac{1}{2}v_{2,2,3}+\frac{1}{2}v_{2,2,5}+\frac{1}{2}v_{2,2,6}+v_{2,3,3}+v_{2,3,5}+\frac{1}{2}v_{3,3,3}-3 .
           \end{eqnarray*}
Define $\tilde{a_{2}}$ to be the last expression in the inequality above.
The above argument can be split into the following three cases.
\begin{itemize}
 \item First, $P$ has at least two cusps.
 \item Second, $P$ has only one cusp and $e_5\geq{2}$. 
 \item Finally, $P$ has only one cusp and $e_5=1$.
\end{itemize}

Let us first consider the case $P$ has at least two cusps. 
By using the identity (9), together with the fact that $P$ has at least one $\frac{\pi}{5}$-edge, we obtain
\[
v_{2,2,5}+v_{2,3,5}\geq{2} .
\]
This inequality means that
\begin{eqnarray*}
\tilde{a_2} &=& (v_{2,2,2,2}+v_{2,3,6}+v_{3,3,3})+\frac{1}{2}v_{2,2,3}+\frac{1}{2}v_{2,2,5}+\frac{1}{2}v_{2,2,6}+v_{2,3,3}+v_{2,3,5}+\frac{1}{2}v_{3,3,3}-3 \\
                   &\geq& 2+\frac{1}{2}v_{2,2,3}+\frac{1}{2}v_{2,2,6}+v_{2,3,3}+1+\frac{1}{2}v_{3,3,3}-3\geq{0} .
\end{eqnarray*}

Next we consider the case that $P$ has only one cusp and $e_5\geq{2}$. 
By using the identity (9), we obtain 
\[
v_{2,2,5}+v_{2,3,5}\geq{4} .
\]
Hence, we obtain the following inequalities.
\begin{eqnarray*}
\tilde{a_2} &=& (v_{2,2,2,2}+v_{2,3,6}+v_{3,3,3})+\frac{1}{2}v_{2,2,3}+\frac{1}{2}v_{2,2,5}+\frac{1}{2}v_{2,2,6}+v_{2,3,3}+v_{2,3,5}+\frac{1}{2}v_{3,3,3}-3 \\
                   &\geq& 1+\frac{1}{2}v_{2,2,3}+\frac{1}{2}v_{2,2,6}+v_{2,3,3}+2+\frac{1}{2}v_{3,3,3}-3\geq{0} .
\end{eqnarray*}

Finally, we consider the case that $P$ has only one cusp and $e_5=1$. 
In this case, $e_5=1$ implies that $(v_{2,2,5}, v_{2,3,5})=(0, 2), (1,1), (2,0)$. 
If $(v_{2,2,5}, v_{2,3,5})=(0, 2)$, we obtain 
\begin{eqnarray*}
\tilde{a_2} &=& (v_{2,2,2,2}+v_{2,3,6}+v_{3,3,3})+\frac{1}{2}v_{2,2,3}+\frac{1}{2}v_{2,2,5}+\frac{1}{2}v_{2,2,6}+v_{2,3,3}+v_{2,3,5}+\frac{1}{2}v_{3,3,3}-3 \\
                   &=& 1+\frac{1}{2}v_{2,2,3}+\frac{1}{2}v_{2,2,6}+v_{2,3,3}+2+\frac{1}{2}v_{3,3,3}-3\geq{0} .
\end{eqnarray*}

If $(v_{2,2,5}, v_{2,3,5})=(1,1)$, we obtain
\begin{eqnarray*}
\tilde{a_2} &=& (v_{2,2,2,2}+v_{2,3,6}+v_{3,3,3})+\frac{1}{2}v_{2,2,3}+\frac{1}{2}v_{2,2,5}+\frac{1}{2}v_{2,2,6}+v_{2,3,3}+v_{2,3,5}+\frac{1}{2}v_{3,3,3}-3 \\
                   &=& 1+\frac{1}{2}v_{2,2,3}+\frac{1}{2}+\frac{1}{2}v_{2,2,6}+v_{2,3,3}+1+\frac{1}{2}v_{3,3,3}-3\geq{-\frac{1}{2}} .
\end{eqnarray*}
Since $\tilde{a}_2$ is an integer, we conclude $\tilde{a}_2\geq{0}$.

Finally, we consider the case $(v_{2,2,5}, v_{2,3,5})=(2,0)$, that is, 
\[
\tilde{a_{2}}={(v_{2,2,2,2}+v_{2,3,6}+v_{3,3,3})+\frac{1}{2}v_{2,2,3}+1+\frac{1}{2}v_{2,2,6}+v_{2,3,3}+\frac{1}{2}v_{3,3,3}-3} .
\]

Here, $P$ has an only one cusp. Thus we can assume that $(v_{2,2,2,2}, v_{2,3,6}, v_{3,3,3})=(1,0,0) \  \text{or} \  (0, 1, 0) \  \text{or} \  (0, 0, 1)$.

If $(v_{2,2,2,2}, v_{2,3,6}, v_{3,3,3})=(0, 1, 0)$, $P$ has a $\frac{\pi}{3}$-edge shared by vertices of type $(2,3,6)$ and $(3, q, r)$. Thus, we obtain
\[
v_{2,2,3}+2v_{2,3,3}\geq{1} .
\]
Furthermore, $P$ also has a $\frac{\pi}{6}$-edge shared by vertices of type $(2,3,6)$ and $(2,2,6)$, so that we obtain
\[
v_{2,2,6}\geq{1} .
\]
These inequalities imply that
\begin{eqnarray*}
\tilde{a_2} &=&(v_{2,2,2,2}+v_{2,3,6}+v_{3,3,3})+\frac{1}{2}v_{2,2,3}+1+\frac{1}{2}v_{2,2,6}+v_{2,3,3}+\frac{1}{2}v_{3,3,3}-3\\
                   &\geq& 1+\frac{1}{2}+1+\frac{1}{2}-3=0 .
\end{eqnarray*}

If $(v_{2,2,2,2}, v_{2,3,6}, v_{3,3,3})=(0,0,1)$, $P$ has three $\frac{\pi}{3}$-edges shared by vertices of type $(3,3,3)$ and $(3,q,r)$. Thus, we obtain
\[
v_{2,2,3}+v_{2,3,3}\geq{3} .
\]
This inequality implies that 
\begin{eqnarray*}
\tilde{a_2} &=&(v_{2,2,2,2}+v_{2,3,6}+v_{3,3,3})+\frac{1}{2}v_{2,2,3}+1+\frac{1}{2}v_{2,2,6}+v_{2,3,3}+\frac{1}{2}v_{3,3,3}-3\\
                   &\geq& 1+\frac{3}{2}+1+\frac{1}{2}v_{2,2,6}+\frac{1}{2}-3>{0} .
\end{eqnarray*}

If $(v_{2,2,2,2}, v_{2,3,6}, v_{3,3,3})=(1,0,0)$, then
\[
\tilde{a_{2}}=\frac{1}{2}v_{2,2,3}+\frac{1}{2}v_{2,2,6}+v_{2,3,3}-1 .
\]

In this case, if $P$ has at least one $\frac{\pi}{3}$-edge or $\frac{\pi}{6}$-edge, we see that $\tilde{a_2}\geq{0}$ .
Therefore, by the assumption, we obtain $\tilde{a_2}\geq{0}$. 
Thus, we can apply Proposition 1 to $H_{2,3,5,6}(t)$ and conclude that the growth rate of $P$ is a Perron number. \qed

If $P$ has no $\frac{\pi}{3}$-edges and $\frac{\pi}{6}$-edges, all the dihedral angles of $P$ are $\frac{\pi}{2}$ and  $\frac{\pi}{5}$.
For this reason, we assume that all the dihedral angles of $P$ are $\frac{\pi}{2}, \frac{\pi}{5}$ and $e_5=1$. 
Then, by means of Steinberg's formula, we can calculate the growth function of $P$ as 
\[
\dfrac{1}{f_{P}(t)}=\dfrac{(t-1)}{[2,2,2,5]}H_{2,5}(t)
\]
where
\begin{eqnarray*}
H_{2,5}(t) &=& (v_{2,2,2,2}-1)t^6+(v_{2,2,2,2}+f-5)t^5+(v_{2,2,2,2}+f-5)t^4 \\
                  &+& (v_{2,2,2,2}+f-5)t^3+(v_{2,2,2,2}+f-5)t^2+(f-5)t-1.
\end{eqnarray*}
This implies that all the coefficients of $H_{2,5}(t)$ except its constant term are nonnegative and the coefficients of $k$-th term for $k=2,3,4,5$ are positive. 
Therefore, we obtain the following proposition. 

\begin{prop}
Suppose $P$ has at least one $\frac{\pi}{5}$-edge. Then, the growth rate of $P$ is a Perron number.
\end{prop}

\subsection{The growth rates of non-compact Coxeter polyhedra whose dihedral angles are $\frac{\pi}{2},\frac{\pi}{3},\frac{\pi}{4}, \frac{\pi}{5}$ and $\frac{\pi}{6}$}
Let $P$ be a non-compact hyperbolic Coxeter polyhedron whose dihedral angels are $\frac{\pi}{2},\frac{\pi}{3}, \frac{\pi}{4},\frac{\pi}{5}$ and $\frac{\pi}{6}$.
We also assume that $P$ has at least one $\frac{\pi}{4}$-edge.
In a manner similar to the above arguments, the growth function of $P$ can be written as

\[
\dfrac{1}{f_{P}(t)}=\dfrac{(t-1)}{[2,4,6,10]}H_{2,3,4,5,6}(t)
\]
where $H_{2,3,4,5,6}(t)$ satisfies the following identity.
\begin{eqnarray*}
H_{2,3,4,5,6}(t)-H_{2,3,5,6}(t) &=& v_{2,2,4}t^{17}+(\frac{1}{2}v_{2,2,4}+\frac{1}{2}v_{2,3,4}+\frac{5}{2}v_{2,4,4})t^{16}\\
                                                      &+& (v_{2,2,4}+\frac{3}{2}v_{2,3,4}+\frac{9}{2}v_{2,4,4})t^{15}+(2v_{2,2,4}+3v_{2,3,4}+\frac{13}{2}v_{2,4,4})t^{14}\\
                                                      &+&(\frac{5}{2}v_{2,2,4}+\frac{9}{2}v_{2,3,4}+\frac{17}{2}v_{2,4,4})t^{13}+(\frac{7}{2}v_{2,2,4}+6v_{2,3,4}+\frac{21}{2}v_{2,4,4})t^{12}\\
                                                      &+&(4v_{2,2,4}+7v_{2,3,4}+\frac{23}{2}v_{2,4,4})t^{11}+(\frac{9}{2}v_{2,4,4}+\frac{15}{2}v_{2,3,4}+12v_{2,4,4})t^{10}\\
                                                      &+&(\frac{9}{2}v_{2,2,4}+\frac{15}{2}v_{2,3,4}+12v_{2,4,4})t^9+(\frac{9}{2}v_{2,4,4}+\frac{15}{2}v_{2,3,4}+12v_{2,4,4})t^8\\
                                                      &+&(\frac{9}{2}v_{2,2,4}+\frac{15}{2}v_{2,3,4}+11v_{2,4,4})t^7+(4v_{2,2,4}+7v_{2,3,4}+\frac{19}{2}v_{2,4,4})t^6\\
                                                      \end{eqnarray*}
\begin{eqnarray*}
                                                      &+&(\frac{7}{2}v_{2,4,4}+6v_{2,3,4}+\frac{15}{2}v_{2,4,4})t^5+(\frac{5}{2}v_{2,2,4}+\frac{9}{2}v_{2,3,4}+\frac{11}{2}v_{2,4,4})t^4\\
                                                      &+&(2v_{2,2,4}+3v_{2,3,4}+\frac{7}{2}v_{2,4,4})t^3+(v_{2,2,4}+\frac{3}{2}v_{2,3,4}+\frac{3}{2}v_{2,4,4})t^2\\
                                                      &+&(\frac{1}{2}v_{2,2,4}+\frac{1}{2}v_{2,3,4}+\frac{1}{2}v_{2,4,4})t-1.
\end{eqnarray*}

The above identity means that if the coefficient of $k$-th term of $H_{2,3,5,6}(t)$ is nonnegative (resp. positive), then the coefficient of $k$-th term of $H_{2,3,4,5,6}(t)$ is nonnegative (resp. positive). 
By using Proposition 1, we obtain the following theorem.

\begin{theo}
Let $P$ be a non-compact hyperbolic Coxeter polyhedron whose dihedral angles are $\frac{\pi}{2}, \frac{\pi}{3}, \frac{\pi}{4}, \frac{\pi}{5}$ and $\frac{\pi}{6}$.
Then the growth rate of $P$ is a Perron number.
\end{theo}

\section{Acknowledgement}
The author is greatly indebted to his supervisor Professor Yohei Komori for his guidance and encouragement in the course of this work.
He also thanks the referee for thorough reading of this paper and for helpful comments and suggestions
which make this paper more readable.


\end{document}